# ATTITUDE QUATERNION ESTIMATION USING A SPECTRAL PERTURBATION APPROACH

## Adam L. Bruce[*]

All quaternion methods for static attitude determination currently rely on either the spectral decomposition of a 4 × 4 matrix (*q*-Method) or finding the maximum eigenvalue of a 4th-order characteristic equation (QUEST). Using a spectral perturbation approach, we show it is possible to analytically estimate the attitude quaternion to high accuracy and recursively calculate the maximum-likelihood attitude quaternion to arbitrary numerical precision. Analytic or recursive estimation removes several numerical difficulties which are inherent to state-of-the-art algorithms, suggesting our results may substantially benefit high frequency and limited memory embedded software implementations, such as those commonly used on spacecraft computers.

## INTRODUCTION

A key component of spacecraft control systems is the algorithm which estimates instantaneous spacecraft attitude. Given a set of body frame vector measurements $\{\mathbf{b}_n\}_n$ with (known) reference frame components $\{\mathbf{r}_n\}_n$, we seek to determine the attitude matrix $A$, or quaternion $q$, best rotating the reference vectors into the body frame. The preeminent framework for attitude determination is the minimization of the Wabha loss function[1]

$$J(A) = \sum_{n=1}^{N} \frac{1}{2} w_n ||\mathbf{b}_n - A\mathbf{r}_n||^2, \qquad (1)$$

where $N$ is the total number of measurement vectors and $w_n$ are weights[*]. Although Wabha's problem appears similar to traditional least-squares minimization, the normalization of measurement vectors means the sample space is topologically equivalent to $S^2$ (the unit sphere). Since $S^2$ is not a Hilbert space, or even a linear space, it need not obey the projection theorem, and the normal equations of traditional least-squares minimization do not exist. Furthermore, minimizing $J(A)$ over all possible attitude matrices is equivalent to functional optimization on $SO_3$, the Lie group of three dimensional rotations. The nontrivial topology of $SO_3$ precludes most simple optimization proceedures.

---

[*]This research was completed while the author was a graduate student at Purdue University and is being released as a university work product qualifying as fundamental research under ITAR 120.11(a)(8). The author is currenty Systems Engineer II, Guidance, Navigation, and Control Center, Raytheon Missile Systems, 1151 E. Hermans Rd, Tucson, AZ 85756.

[*]In a slight abuse of terminology, we shall refer to the set $\{\mathbf{b}_n\}_n^N \cup \{\mathbf{r}_n\}_n^N$ as the *vector measurements* of the attitude determination algorithm, despite the fact that only the $\mathbf{b}_n$ are measured in real-time, with the $\mathbf{r}_n$ typically calculated from a model.

Despite being more complex, there exist several efficient algorithms to solve Wabha's problem. For all $N \geq 2$ the solution is unique, as there are more degrees of freedom ($2n$) than independent unknown parameters (3), and the loss function is convex. Of these algorithms, *matrix methods* solve for the attitude matrix directly and *quaternion methods* solve for the equivalent attitude quaternion. The latter class can only resolve the quaternion up to $\pm 1$, since $A(q) = A(-q)$.

A classic result of P. Davenport (as reported by Wertz[2] and more recently by Markley and Crassidis[3]) is that solving Wabha's problem is equivalent to solving

$$Kq = \lambda_m q, \qquad (2)$$

where $q$ is the (Wabha) optimal quaternion, $K$ is a $4 \times 4$, symmetric, traceless matrix, whose entries are functions of the weights and vector measurements exclusively, and $\lambda_m$ is the largest eigenvalue of $K$. An estimate for $\lambda_m$ can be used to easily calculate the optimal quaternion, so an efficient method to approximate $\lambda_m$ completes the solution of Wabha's problem.

In this work we present a method for approximating $\lambda_m$ using a perturbation approach. The first-order truncation of the perturbation solution produces an approximation for $\lambda_m$ which is both simple and accurate to roughly $\mathcal{O}(10^{-3})$, and can also be extended into a fast-converging recursive formula for higher accuracies. Higher-order truncations of the perturbative series are also investigated.

## QUATERNION ATTITUDE ESTIMATION

We assume the reader is familiar with using quaternions for attitude analysis. In this section we establish our notation and conventions and review the concepts necessary for understanding the rest of the paper.

### Attitude Quaternions

Quaternions are elements of the algebra of hypercomplex numbers generated by the symbols $\{1, i, j, k\}$ and multiplication table

$$i^2 = j^2 = k^2 = ijk = -1. \qquad (3)$$

The quaternion algebra is denoted $H$ after their inventor, W. R. Hamilton.[4] For any $q \in H$ there are always four reals $r_0, r_1, r_2$ and $r_3$ for which

$$q = r_0 + r_1 i + r_2 j + r_3 k. \qquad (4)$$

$r_0$ is called the *scalar part* of $q$, while $r_1 i + r_2 j + r_3 k$ is called the *vector part*. The set $P = \{q \in H : r_0 = 0\}$ is called the set of *pure quaternions*.

We define projectors $S : H \to R$ and $V : H \to P$ where

$$Sq = r_0, \quad Vq = r_1 i + r_2 j + r_3 k \quad \text{for all} \quad q \in H. \tag{5}$$

Any quaternion can be written $q = Sq + Vq$. We denote $Vq_a = r_a$, $a = 1, 2, 3$. Any binary operation $\star$ of quaternions $p$ and $q$ is uniquely defined by $S(p \star q)$ and $V(p \star q)$. Let $p, q \in H$. Then

$$S(p + q) := Sp + Sq, \quad V(p + q) := \sum_a (Vp_a + Vq_a). \tag{6}$$

Using Eq. (3), one can check $qp$ is equivalent to

$$S(qp) = SqSp - Vq \cdot Vp, \quad V(qp) = SqVp + SpVq + Vq \times Vp, \tag{7}$$

where $Vq \cdot Vp = \sum_a Vq_a Vp_a$ and $Vq \times Vp = \sum_{b,c} \epsilon_{abc} Vq_b Vp_c$. The conjugate, $q^*$, of $q$ is $q^* = Sq - Vq$. As with complex numbers the norm is $|q| = \sqrt{q^* q}$.

The *unit quaternions* are the multiplicative subgroup $H_U := \{q \in H : |q| = 1\}$. We find $|q| = 1$ is equivalent to $(Sq)^2 + \sum_a (Vq_a)^2 = 1$. Although the use of quaternions to represent attitude is well-known in the aerospace community, the reason why this is possible is less so. This is shown by the following

**Theorem.** $H_U \simeq SU_2$.

*Proof.* Let $\varphi : H \to SU_2$ send $1, i, j,$ and $k$ to $I, i\sigma_1/2, i\sigma_2/2,$ and $i\sigma_3/2$, where the $\sigma_a$ are the Pauli matrices. $\varphi$ is a homomorphism because $I, i\sigma_1/2, i\sigma_2/2,$ and $i\sigma_3/2$ preserve the multplication table of the quaternion basis. Furthermore, the condition that $|q| = 1$ is necessary and sufficient for $\varphi(q) \in SU_2$, establishing the bijection. □

The existence of the covering map $K : SU_2 \to SO_3$, obtained under the conjucation of the quaternion roots of $-1$, is similarly well known[*], and the correspondence between $i, j, k$ and the unit basis of $R^3$ shows this is equivalent to proper rotation of vectrices, hence a parameterization of $SO_3$. That is to say $(K \circ \varphi)q \equiv A(q) \in SO_3$. Furthermore ker $K = \{\pm 1\}$ so $A(q) = A(-q)$, the well known fact that rotations corresponding to a quaternion $q$ and $-q$ are identical. As a universal covering map, $K$ must be a homomorphism, so $A = K \circ \varphi$ is also a homomorphism, which is to say $A(pq) = A(q)A(p)$ for all unit quaternions $q$ and $p$. This fact implies unit quaternion multiplication is equivalent to composing the rotations corresponding to each quaternion in the product. By induction it follows that $q_1 q_2 q_3 \cdots q_n$ corresponds to $A(q_1)A(q_2) \cdots A(q_n)$, for unit quaternions $q_1, q_2, q_3, \cdots, q_n$.

---

[*]This is commonly called the "quaternion attitude matrix" or a similar name in the engineering literature. To obtain it send $i, j, k$ to $qiq^*, qjq^*, qkq^*$.

Finally we show that any quaternion admits the representation $re^p$, where $r$ is real, $p$ is pure, and the exponential map is defined as

$$e^p = \sum_{a=0}^{\infty} \frac{p^a}{a!}, \quad p^0 \equiv 1. \tag{8}$$

**Theorem.** *If $q \in H$ then there are $r \in R$ and $p \in P$ for which $q = re^p$.*

*Proof.* $H = R \times H_U$ because $q = |q|(q/|q|)$. Let $\text{Lie}: G \to \mathfrak{g}$ be the map which sends a Lie group to its algebra. Since $H_U \simeq SU_2$ and $SU_2$ is compact, the claim is proven if $\text{Lie}(H_U) = P$. Lie algebras of isomorphic Lie groups are themselves isomorphic, so we need only show $\mathfrak{su}_2 \simeq P$. Define $\varphi|_P : P \to \mathfrak{su}_2$ as the map which sends $i, j$ and $k$ to $i\sigma_1/2$, $i\sigma_2/2$, and $i\sigma_3/2$ respectively. Since $[i\sigma_a/2, i\sigma_b/2] = i\epsilon_{abc}\sigma_c/2$, the Pauli matrices generate $\mathfrak{su}_2$. Furthermore, the algebra $[\varphi|_P^{-1}(i\sigma_a/2), \varphi|_P^{-1}(i\sigma_b/2)] = \epsilon_{abc}\varphi|_P^{-1}(i\sigma_c/2)$ is clearly $P$. So $\varphi|_P : P \to \mathfrak{su}_2$ is an isomorphism. □

**Quaternion Perturbation Expansions**

Only unit quaternions can be used to represent attitude. We therefore need a way of calculating the expansion $pq = q_0 + q_1 + q_2 + \cdots$ when $p, q \in H_U$.

Suppose $q, q_0 \in H_U$. There is always a unique $Q \in H_U$ for which $q = Qq_0$, namely $Q = qq_0^*$. But since $Q \in H_U$ there is a $p_Q \in P$ for which $Q = e^{p_Q}$. Then

$$q = e^{p_Q} q_0 = \left(1 + p_Q + \frac{1}{2}p_Q^2 + \cdots\right) q_0 = q_0 + p_Q q_0 + \frac{1}{2}p_Q^2 q_0 + \cdots \tag{9}$$

so $q_a = p_Q^a q_0 / a!$. Note that the $a^{\text{th}}$ term in the perturbation series is $\mathcal{O}(|p_Q|^a)$. If $|p_Q| \ll 1$, viz. $Q$ is small, then the approximation

$$q \approx q_0 + p_Q q_0 \tag{10}$$

is good. As this truncates Eq. (9) to $a = 1$ we refer to Eq. (10) as the *first-order* expansion of Eq. (9).

**Vector Treatment of Quaternions**

Quaternions are often referred to as "four-dimensional vectors" in aerospace texts. Strickly speaking this is false*, but there are situations when it is useful represent to quaternions as vectors and use linear algebra. We represent $q \in H$ as a vector $q \in R^4$ by sending the quaternion basis to the Euclidean basis. In a slight abuse of notation we define the quaternion vector as

$$q = \begin{bmatrix} Vq_1 \\ Vq_2 \\ Vq_3 \\ Sq \end{bmatrix} \equiv \begin{bmatrix} Vq \\ Sq \end{bmatrix}, \tag{11}$$

---

*The algebraic structure fails to be preserved, despite the fact that $H$ is homeomorphic to $R^4$ and $H_U$ to $S^3$.

in accordance with the typical aerospace convention. When linear algebra operations are performed on quaternions, the reader should assume they are represented as vectors.

**Davenport's Eigenproblem**

Let us now return to Wabha's loss function. For normalized vector measurements

$$J(A) = \sum_n \frac{1}{2} w_n (\mathbf{b}_n - A\mathbf{r}_n)^T (\mathbf{b}_n - A\mathbf{r}_n) = \lambda_0 - G(A), \tag{12}$$

where $\lambda_0 = \sum_n w_n$ and $G(A) = \sum_n w_n \mathbf{b}_n^T A \mathbf{r}_n$. Davenport showed $G(A) = q^T K q$, where

$$K = \begin{bmatrix} \rho - \sigma I & \mathbf{z} \\ \mathbf{z}^T & \sigma \end{bmatrix}, \tag{13}$$

and

$$B = \sum_{n=1}^{N} w_n \mathbf{b}_n \mathbf{r}_n^T, \quad \rho = B + B^T, \quad \mathbf{z} = \sum_{n=1}^{N} w_n \mathbf{b}_n \times \mathbf{r}_n, \text{ and } \sigma = \text{tr} B. \tag{14}$$

Markley and Crassidis[3] derive this result. Recall Wabha's problem is convex, so minimizing $J$ is equivalent to maximizing $G$. Since $K$ is symmetric, the spectral theorem guarantees the existence of an orthonormal basis of $R^4$ consisting only of eigenvectors of $K$. This means

$$q^T K q = \sum_{\lambda \in \Lambda(K)} |a_\lambda|^2 \lambda \leq \lambda_m, \tag{15}$$

where $\Lambda(K)$ is the spectrum of $K$. One finds the minimum of $J$ is $\lambda_0 - \lambda_m$ and the quaternion minimizing $J$ obeys $Kq = \lambda_m q$.

**The QUEST Algorithm**

Davenport's eigenproblem can be solved by directly computing the spectral factorization of $K$, which is called the *Davenport q-Method* or simply the *q-Method*. The Spectral factorization of a $4 \times 4$ matrix can be too computationally expensive to be used on an embedded platform. Markley and Crassidis[3] record how and this limitation became apparent during the 1970's, when the High Energy Astronomy Observatory (HEAO-B) needed faster attitude updates than the *q*-method could deliver. Though computing power has increased since then, limitations are still present today, such as on small spacecraft.

Shuster proposed the QUEST algorithm[5,6,7] as a faster alternative to the *q*-Method. Since the Rodrigues parameters, defined as $p = q/Sq$, recover the quaternion via

$$q = \frac{1}{\sqrt{1 + p^T p}} \begin{bmatrix} p \\ 1 \end{bmatrix}, \tag{16}$$

we find $q$ by finding $p$. The Davenport eigenproblem is rewritten as

$$[\rho - (\sigma + \lambda_m)I]p = -\mathbf{z}, \quad \mathbf{z}^T p + \sigma = \lambda_m, \tag{17}$$

so

$$p = -[\rho - (\sigma + \lambda_m)I]^{-1}\mathbf{z}, \tag{18}$$

with $\lambda_m$ obtained by solving

$$\mathbf{z}^T[\rho - (\sigma + \lambda_m)I]^{-1}\mathbf{z} + \lambda_m - \sigma = 0. \tag{19}$$

The $3 \times 3$ matrix inverse is calculated by an analytic formula, so the only iteration involved is the nonlinear root finding. $J \approx 0$ when minimized, from which it follows that $\lambda_m \approx \lambda_0$. This justifies choosing $\lambda_0$ to initialize the nonlinear root finder to solve Eq. (19).

Directly approximating $\lambda_m$ by $\lambda_0$ is a zeroth-order estimate, which is sometimes good enough. Since it doesn't incorporate measurement knowelege however, it is generally too inprecise when measuring instruments are good. Another reason to find a more precise estimate is to calculate Shuster's $\chi^2$-distributed TASTE statistic,[8]

$$\text{TASTE} = \lambda_0 - \lambda_m, \tag{20}$$

which can be used for real-time data validation.

The following section derives a first-order method for approximating $\lambda_m$. This gives approximations for both $q$ and TASTE which can be calulated without iteration and using only floating-point operations.

### PERTURBATIVE SOLUTION OF THE DAVENPORT EIGENPROB­LEM

The main idea is to use Eq. (10) to solve the Davenport Eigenproblem by a perturbation series. To do this we also need peturbation expansions for $\lambda_m$ and $K$. The former is just the expansion $\lambda_m = \lambda^{(0)} + \lambda^{(1)} + \lambda^{(2)} + \cdots$ of a scalar, while we denote the later $K = K_0 + \delta_K$.

Now suppose we identify $\lambda^{(0)}$ with $\lambda_0$ and drop all $\mathcal{O}(\geq 2)$ terms. This sets the base rotation as the quaternion corresponding to the estimate $\lambda_m \approx \lambda_0$, which we seek to refine by a smaller rotation. This also allows $q$ to be written as $q_0 + pq_0$, where $p$ is the pure generator of the perturbing rotation and $q_0$ is defined as is the quaternion with Rodrigues parameters

$$p_0 = -[\rho - (\sigma + \lambda_0)I]^{-1}\mathbf{z}. \tag{21}$$

We define $K_0$ as the "synthetic" $K$-matrix for which, given $\lambda_0$ and $q_0$

$$K_0 q_0 = \lambda_0 q_0. \tag{22}$$

That is, if $\lambda_m = \lambda_0$ and $q = q_0$ exactly, the $K$ constructed by measurements would be $K_0$. Note that this means $K_0$ must be symmetric and traceless. Since we assume $\lambda_0$ and $q_0$ are close to but not exactly $\lambda_m$ and $q$, we define the matrix $\delta_K = K - K_0$, whose entries are first-order perturbation quantities. Now,

$$(K_0 + \delta_K)(q_0 + pq_0) = (\lambda_0 + \lambda^{(1)})(q_0 + pq_0). \tag{23}$$

The left side is

$$K_0 q_0 + \delta_K q_0 + K_0 p q_0 + \delta_K p q_0 \approx K_0 q_0 + \delta_K q_0 + K_0 p q_0, \tag{24}$$

since $\delta_K p$ is second-order, while the right is

$$\lambda_0 q_0 + \lambda^{(1)} q_0 + \lambda_0 p q_0 + \lambda^{(1)} p q_0 \approx \lambda_0 q_0 + \lambda^{(1)} q_0 + \lambda_0 p q_0. \tag{25}$$

Removing the zeroth-order terms by enforcing Eq. (22),

$$\delta_K q_0 + K_0 p q_0 = \lambda^{(1)} q_0 + \lambda_0 p q_0. \tag{26}$$

And since $q_0^T K_0 = (q_0 K_0)^T = \lambda_0 K_0$, we premultiply by $q_0^T$ and find

$$\lambda^{(1)} = q_0^T \delta_K q_0 = q_0^T (K - K_0) q_0 = q_0^T K q_0 - \lambda_0. \tag{27}$$

So

$$\lambda_m \approx \lambda_0 + \lambda^{(1)} = q_0^T K q_0. \tag{28}$$

The first-order estimate of the optimal quaternion, which we denote $q_1$, therefore has Rodrigues paramters $[(\sigma + q_0^T K q_0)I - \rho]^{-1}\mathbf{z}$.

**First-Order Performance**

We analyze the first-order performance using Monte Carlo simulations for the maximum eigenvalue, similar to those previously done by Bruce and Frueh (2016).[9] In this work, the error distribution parameter for the Monte Carlo simulations was the modulus of the error rotation vector. In this case we are able to report distributions in terms of the more conventional error parameter

$$\lambda_1 - \lambda_m. \tag{29}$$

Monte Carlo simulations were performed with $5 \times 10^6$ points. The histogram binning method used was identical to the previous work with the histogram density parameter was again selected to be 3333 for a good compromise between manageable runtime and good distributional reconstruction. We have inluded Fig. 1 for readers unfamiliar with the previous work to illustrate the influence of the histogram density parameter $\rho_H$* on the quality of the histogram reconstruction using a Gaussian distribution with $\mu = 0$ and $\sigma = 0.1$.

---

*Not to be confused with the matrix $\rho = B + B^T$.

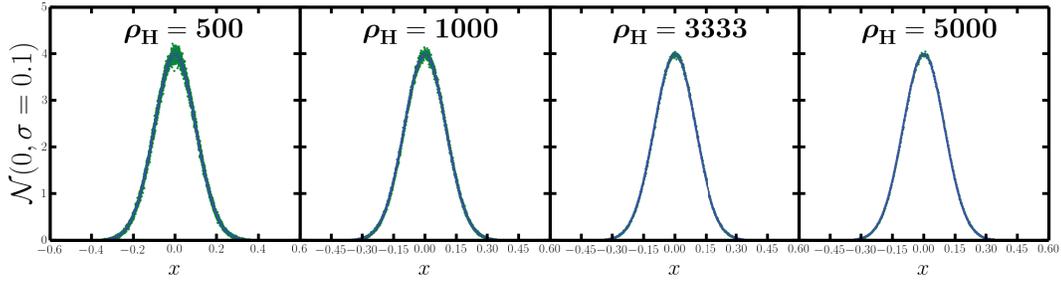

**Figure 1.** Convergence of histogram reconstruction (green) to the underlying population distribution (blue) as $\rho_H$ is increased. $N = 5 \times 10^6$ samples of $\mathcal{N}(0, 0.1)$ with $N_{\mathbf{Bins}} = 10^4,\ 5\times 10^3,\ 1.5\times 10^3,$ and $10^3$ respectively (*reproduced from Bruce and Frueh.[9]*)

As in the previous work, we performed simulations with angular noise introduced to two vector measurements, and allowed the $\sigma$'s of the noise distributions to take on different values, so that the pair $(\sigma_1, \sigma_2)$ could identify the simulation. The pairs chosen for this study were

$$(0.1°, 0.1°), (0.1°, 0.5°), (0.1°, 1.0°), (0.5°, 0.5°), (0.5°, 1.0°), \text{ and } (1.0°, 1.0°). \quad (30)$$

The distributions for $\lambda_m - \lambda_1$ are shown in Fig. 2.

A first feature we see is the distributions are exponential. This is expected as $\lambda_1 \leq \lambda_m$ by construction. Second, we find the probability mass of the distribution is overwhelmingly within error parameters of roughly $10^{-3}$ or less. Those cases in which the error parameter is larger result from poor measurements, which violate the assumption that the true eigenvalue is close to $\lambda_0$ and would be thrown out using TASTE. In some cases the bin separation is visible. This is not due to coarseness of the histogram density, but rather because the probability mass contained by small error bins is so much larger than larger error parameter bins containing only outliers. To obtain a smooth distribution these outlying values would need to be discarded using TASTE. One finds in all cases that the probability mass of the smallest error bin is much larger than all other bins. These simulations show for a variety of cases that the approximation of $\lambda_m = \lambda_1$ is good to at least $10^{-3}$ (excepting outlying cases) and generally more.

**Higher-Order Perturbations**

The second-order expansion of the perturbation series is $q = q_0 + pq_0 + p^2 q_0/2$ in the quaternion and $\lambda = \lambda_0 + \lambda^{(1)} + \lambda^{(2)}$ in the eigenvalue. A similar proceedure as the first-order case leads to an equation in which $p$ is still present. An approximate solution is retrieved by incorporating the first-order quaternion, $q_1 = q_0 + pq_0$, so that

$$p \approx q_1 q_0^* - 1. \quad (31)$$

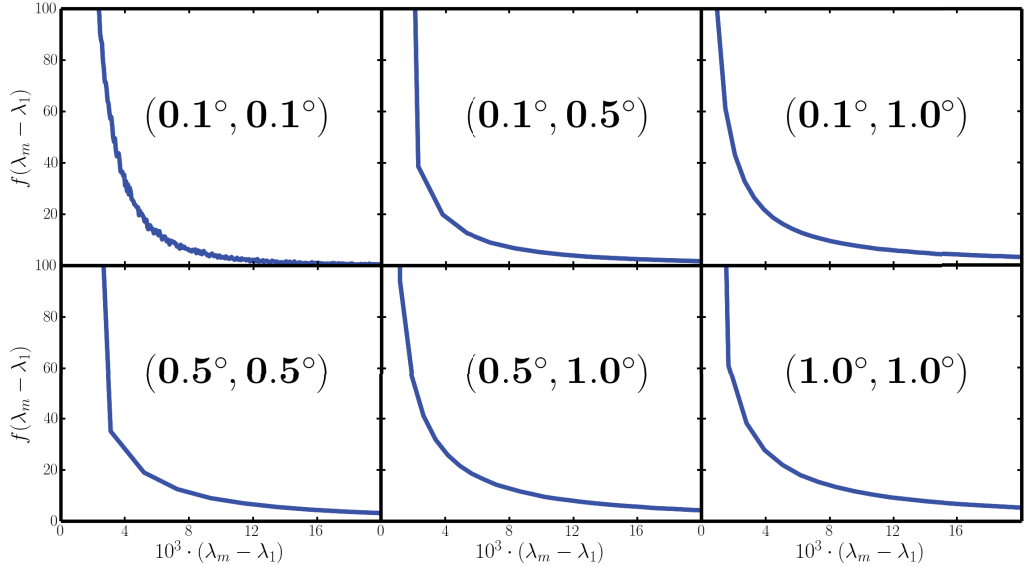

**Figure 2.** Monte Carlo histograms of the six simulations. The error parameter is exponentially distributed with the vast majority of probability mass accounted for by errors of roughly $10^{-3}$ or less.

Although this leads to a solution for the eigenvalue, it is much more complicated than the first-order expansion and we defer the study of these expansions to a subsequent work.

## RECURSIVE ALGORITHM

The calculation of $\lambda_1$ by perturbation methods does not solve issues such as the need to find the inverse of $(\lambda + \sigma) - \rho$ or estimating $\lambda$ to higher precision. The latter of these issues may be solved by formulaing the recursive process

$$\lambda_{a+1} = q_a^T K q_a, \qquad (32)$$

since each $\lambda_a$ may be regarded as a closer estimate of $\lambda$ than the last, the perturbation assumption gets better with each iteration, as $J_a < J_{a-1}$ for each $a$ and the quaterion rotating each $q_a$ to its successor $q_{a+1}$ is progressively smaller. So long as the initial perturbation is satisfied (equivalent to the vector measurements having good SNR), this process will converge to the true solution. In practice 3-4 iterations typically yields convergence to ten or twelve decimal places.

The final difficulty is the sequential inversion of $(\lambda + \sigma) - \rho$. Define $D(\lambda) = [(\lambda + \sigma) - \rho]^{-1}$. We will show it is possible to find a function $f$ for which $D(\lambda_{a+1}) \approx f(\lambda_{a+1} - \lambda_a, D(\lambda_{a+1}))$.

We first state the Neumann or Resolvent series[10, 11] of the operator $A$ (assumed to be on the space $X$),

$$R(A) = (1-A)^{-1} = \sum_{k=0}^{\infty} A^k. \quad (33)$$

This is a generalization of the Taylor series of $1/(1-x)$ to arbitrary linear operators. Let $\|A\|$ be the operator norm on $X$ defined by $\sup_{x \in X} |Ax|/|x|$,[11,10,12] where $|\cdot| : X \to R$ is the vector norm on $X$. $\|\cdot\| : X^X \to R$ is easily shown to be submultiplicative, e.g. $\|AB\| \leq \|A\|\|B\|$. Thus

$$\left\| \sum_{k=0}^{\infty} A^k \right\| \leq \sum_{k=0}^{\infty} \|A^k\| \leq \sum_{k=0}^{\infty} \|A\|^k, \quad (34)$$

whence the condition $\|A\| < 1$ is sufficient for the convergence of $\sum_{k=0}^{\infty} A^k$. Furthermore, this shows the truncation $\sum_{k=0}^{n} A^k$ is good to order $\|A\|^{n+1}$.

Consider the operators $A$ and $B$ for which $\|B\| \ll \|A\| < \|A\| + \|B\| < 1$ and the inverse $(A+B)^{-1}$. We may approximate $(A+B)^{-1}$ knowing $A^{-1}$ by writing $(A+B)^{-1} = (1+A^{-1}B)^{-1}A^{-1}$ and applying the Neumann series to find

$$(A+B)^{-1} = \sum_{k=0}^{\infty} (-A^{-1}B)^k A^{-1} = A^{-1} - A^{-1}BA^{-1} + \mathcal{O}(\|B\|^2). \quad (35)$$

Let us apply the first order result to $D(\lambda_a)$. Since $\lambda_0$ is arbitrary, one may always find weights such that $\|D(\lambda)\| < 1$. Thus

$$D(\lambda_a) = [(\lambda_{a-1} + (\lambda_a - \lambda_{a-1}) + \sigma) - \rho]^{-1} = [D(\lambda_{a-1})^{-1} + (\lambda_a - \lambda_{a-1})]^{-1} \quad (36)$$

or

$$D(\lambda_a) \approx D(\lambda_{a-1}) - (\lambda_a - \lambda_{a-1})D^2(\lambda_{a-1}), \quad (37)$$

which is good to $(\lambda_a - \lambda_{a-1})^2$, a second-order quantity in the perturbation assumption.

This result allows us to write the algorithm in a totally iterative way. Assume that the zeroth-order quantities, $\lambda_0, q_0, D(\lambda_0)$ as well as $K$ and $\mathbf{z}$ are calculated. We then have for $a \geq 1$ the iteration

1. $\lambda_a = q_{a-1}^T K q_{a-1}$,

2. $D(\lambda_a) = D(\lambda_{a-1}) - (\lambda_a - \lambda_{a-1})D^2(\lambda_{a-1})$,

3. $q_a = \left(1 + \mathbf{z}^T D^T(\lambda_a) D(\lambda_a) \mathbf{z}\right)^{-1/2} \begin{bmatrix} D(\lambda_a)\mathbf{z} \\ 1 \end{bmatrix}$,

whence we conclude that (1) only initial quantities need be calculated by inverses, (2) the iteration can be extended indefinitely to produce the true result.

## CONCLUSION

We have presented a method for recursive estimation of the attitude quaternion using a spectral perturbation approach. We verified that even a single iteration of this method, which produces a simple formula for the maximum eigenvalue, yields approximations for the true eigenvalue which are good to the order of $10^{-3}$. We also presented a way to recursively estimate the quaternion from an initial (analytic) inversion of the $3\times 3$ matrix $(\sigma+\lambda_0)I-\rho$. When coupled, these two recursion formulas can be used to implement an attitude determination algorithm which improves on the existing QUEST algorithm by removing the sequential matrix inversions necessary to calculate the maximum eigenvalue and quaternion. Furthermore, simply using the first iteration of the recursion formula provides a nearly analytic estimate of the overall solution. In addition to the applications in attitude determination, the topics of this work also provoke a new avenue of research into perturbation expansions of the quaternions, and their role in spectral approximation for linear operators of the Hamilton algebra.

## ACKNOWLEDGMENTS

The author would like to thank Prof. Carolin Frueh for her comments on this draft and Purdue University for allowing the manuscript to be released.